\documentclass[11pt,a4paper]{article}
\usepackage{ifthen,latexsym,amssymb,amsmath,bbm,fixmath}
\usepackage[nobysame,initials]{amsrefs}
\usepackage[shortlabels]{enumitem} 
\usepackage{graphicx}

\newcommand{\C}[1]{{\protect\mathcal{#1}}}
\newcommand{\B}[1]{{\bf #1}}
\newcommand{\I}[1]{{\mathbbm #1}}

\newcommand{\e}{\varepsilon}

\newcommand{\me}{{\mathrm e}}
\renewcommand{\mid}{:}

\renewcommand{\ge}{\geqslant}
\renewcommand{\le}{\leqslant}

\newcommand{\dd}{\,\mathrm{d}}
\newcommand{\ex}{\mathrm{ex}}

\newif\ifnotesw\noteswtrue

\newcommand{\hide}[1]{}

\newcommand{\beq}[1]{\begin{equation}\label{#1}}
\newcommand{\eeq}{\end{equation}}

\newtheorem{theorem}{Theorem}
\newtheorem{lemma}[theorem]{Lemma}
\newtheorem{proposition}[theorem]{Proposition}
\newtheorem{problem}[theorem]{Problem}

\newcommand{\bpf}[1][Proof.]{\smallskip\noindent{\it #1} }
\newcommand{\qed}{\nolinebreak\mbox{\hspace{5 true pt}%
  \rule[-0.85 true pt]{3.9 true pt}{8.1 true pt}}}
\newcommand{\epf}{\qed \medskip}

\newtheorem{claim}[theorem]{Claim}
\newcommand{\cqed}{\nolinebreak\mbox{\hspace{5 true pt}%
  \rule[-0.85 true pt]{2.0 true pt}{8.1 true pt}}}
\newcommand{\bcpf}{\bpf[Proof of Claim.]}
\newcommand{\ecpf}{\cqed \medskip}

\newcommand{\OPdata}{Oleg Pikhurko\footnote{Supported 
by ERC Advanced Grant 101020255 and Leverhulme Research Project Grant RPG-2018-424.
 }\\
Mathematics Institute and DIMAP\\
University of Warwick\\
Coventry CV4 7AL, UK}

\newcommand{\copex}{\mathrm{co^{+}ex}}
\newcommand{\coex}{\mathrm{co\mbox{-}ex}}
\renewcommand{\r}[1]{r(#1)}
\newcommand{\R}[2]{r(#1,#2)}
\newcommand{\rb}[1]{r[#1]}
\newcommand{\rl}[1]{r_<(#1)}
\newcommand{\rlb}[1]{r_<[#1]}
\newcommand{\Sym}[1]{\mathrm{Sym}_{#1}}
\renewcommand{\deg}{\mathrm{deg}}

\newcommand{\essinf}{\operatorname{ess-inf}}
\newcommand{\Edge}[2]{\mathrm{Edge}_{#1,#2}}

\usepackage{stmaryrd}
\newcommand{\eval}[1]{\left\llbracket#1\right\rrbracket}

\DeclareMathOperator*{\einf}{ess\mbox{-}inf}
\newcommand{\degree}{degree}

\author{\OPdata}
\title{On the limit of the positive $\ell$-\degree\ Tur\'an problem}

\begin{document}

\maketitle

\begin{abstract}
The \emph{minimum positive $\ell$-degree} $\delta^+_{\ell}(G)$ of a non-empty $k$-graph $G$ is
the maximum $m$ such that every $\ell$-subset of $V(G)$ is contained in either none or at least $m$ edges of~$G$;  let $\delta^+_{\ell}(G):=0$ if $G$ has no edges.
For a family $\C F$ of $k$-graphs, let $\copex_\ell(n,\C F)$ be the maximum of $\delta^+_{\ell}(G)$ over all $\C F$-free $k$-graphs $G$ on $n$ vertices. We prove that the ratio $\copex_\ell(n,\C F)/{n-\ell\choose k-\ell}$ tends to limit as $n\to\infty$, answering a question of Halfpap, Lemons and Palmer. Also, we show that the limit can be obtained as the value of a natural 
optimisation problem for $k$-hypergraphons; in fact, we give an alternative description of the set of possible accumulation points of almost extremal $k$-graphs.
\end{abstract}

\section{Introduction}

A \emph{$k$-graph} $G$ is a pair $(V(G),E(G))$, where $V(G)$ is the \emph{vertex set} of $G$ and $E(G)$ is a collection of $k$-subsets of $V(G)$, called \emph{edges}. We call $G$ \emph{non-empty} if $E(G)\not=\emptyset$. 

Fix an integer $\ell$ with $0\le \ell\le k-1$. The \emph{minimum positive $\ell$-\degree} of a non-empty $k$-graph $G$, denoted by $\delta^+_{\ell}(G)$, is
the maximum $m$ such that every $\ell$-subset $L$ of $V(G)$ is contained in either none or at least $m$ edges of~$G$. If $G$ has no edges then we define $\delta^+_{\ell}(G):=0$.
For a family $\C F$ of $k$-graphs, the \emph{positive $\ell$-\degree\ Tur\'an problem} is to determine $\copex_\ell(n,\C F)$, the maximum of $\delta^+_{\ell}(G)$ over all $\C F$-free $k$-graphs $G$ on $n$ vertices. Let 
 \beq{eq:Gamma+}
 \gamma_\ell^+(\C F):=\limsup_{n\to\infty} \frac{\copex_\ell(n,\C F)}{{n-\ell\choose k-\ell}}.
 \eeq

If $\ell=0$, then $\copex_0(n,\C F)$ is  the usual \emph{Tur\'an function} $\ex(n,\C F)$, the maximum number of edges in an $n$-vertex $\C F$-free $k$-graph and thus $\gamma_0^+(\C F)$ is the \emph{Tur\'an density} $\pi(\C F):=\lim_{n\to\infty} \ex(n,\C F)/{n\choose k}$, where the existence of the limit was established by Katona, Nemetz and Simonovits~\cite{KatonaNemetzSimonovits64} by an easy averaging argument. Also, one can show that $\gamma_1^+(\C F)=\pi(\C F)$. (For example, this easily follows from Proposition~\ref{pr:LimExists} which states that the ratio in the right-hand side of~\eqref{eq:Gamma+} tends to a limit.) For surveys of the hypergraph Tur\'an problem, we refer the reader to~Sidorenko~\cite{Sidorenko95} and Keevash~\cite{Keevash11}. 

Another relative is the \emph{$\ell$-\degree\ Tur\'an fuction}
$\coex_\ell(n,\C F)$, the maximum $m$ such that there is an $\C F$-free $n$-vertex $k$-graph $G$ such that for every $\ell$-subset of $V(G)$ we have $\deg_L(G)\ge m$, 
where the \emph{\degree} $\deg_L(G)$ of a set $L\subseteq V(G)$ is the number of edges of $G$ that contain~$L$. Trivially, we always have $\copex_\ell(n,\C F)\ge \coex_\ell(n,\C F)$.
A systematic study of this function for $\ell=k-1$ was started by Mubayi and Zhao~\cite[Proposition~1.2]{MubayiZhao07} who in particular proved that the limit
 \beq{eq:CoexLim}
 \gamma_\ell(\C F):=\lim_{n\to\infty} \frac{\coex_\ell(n,\C F)}{{n-\ell\choose k-\ell}}
 \eeq
 exists for $\ell=k-1$. The existence of the limit in~\eqref{eq:CoexLim} for every $\ell$ was proved by Lo and Markstr\"om~\cite[Proposition~1.5]{LoMarkstrom14} (see also Keevash~\cite[Page 118]{Keevash11}).  Balogh, Clemen and Lidicky~\cite{BaloghClemenLidicky22} present a survey of these (and some related) Tur\'an-type problems for $3$-graphs.

The problem of determining $\copex(n,\C F)$ for $\ell=k-1$ was introduced and studied by
Halfpap, Lemons and Palmer~\cite{HalfpapLemonsPalmer22arxiv}, motivated by an earlier paper of 
Balogh, Lemons and Palmer~\cite{BaloghLemonsPalmer20} who studied positive \degree\ in the context of intersection families. Note that our definition of the minimum positive $\ell$-\degree\ of $G$ deviates from the one in~\cite{BaloghLemonsPalmer20} when $G$ has no edges: in this case the authors of~\cite{BaloghLemonsPalmer20} leave $\delta_\ell^+(G)$ undefined while we set $\delta_\ell^+(G):=0$ (with this leading to slightly cleaner statements of some of our results).

Halfpap, Lemons and Palmer~\cite{HalfpapLemonsPalmer22arxiv} computed $\gamma_2(\{F\})$ (and, in some cases, the exact value of $\copex_{2}(n,\{F\})$) for a few natural $3$-graphs $F$; see also Wu~\cite{Wu22arxiv} for some further such results. The  value of $\gamma_2^+(\{F\})$ is in general different from  $\gamma_2(\{F\})$ (as 
well as from $\pi(\{F\})$).

Halfpap, Lemons and Palmer~\cite{HalfpapLemonsPalmer22arxiv} asked if for every $k$-graph $F$ the ratio $\copex_{k-1}(n,\{F\})/n$ tends to a limit as $n\to\infty$ and proved this (\cite[Proposition~21]{HalfpapLemonsPalmer22arxiv}) in a special case when  every $(k-1)$-set of vertices of $F$ is covered by an edge. 

Here we answer this question for every  $k$-graph family $\C F$ (and an arbitrary integer $\ell$). 

\begin{proposition}\label{pr:LimExists} 
 For every  (possibly infinite) $k$-graph family $\C F$ and every integer $\ell$ with $0\le \ell\le k-1$, the ratio $\copex_\ell(n,\C F)/{n-\ell\choose k-\ell}$ tends to a limit as $n\to\infty$.
\end{proposition}

Our second result, Theorem~\ref{th:main}, describes the set of possible $k$-hypergraphons that are the limits of sequences of almost optimal constructions for $\copex_\ell(n,\C F)$ as $n\to\infty$, in particuar giving
an optimisation problem for $k$-hypergraphons that produces~$\gamma_\ell^+(\C F)$. This can be viewed as the natural limit version of the positive $\ell$-\degree\ Tur\'an problem.

Since the statement of   Theorem~\ref{th:main} is  technical and requires quite a few definitions (for a non-expert), we state it only after presenting the (purely combinatorial) proof of Proposition~\ref{pr:LimExists}.

\section{Proof of Proposition~\ref{pr:LimExists}}\label{LimExists}

The main proof idea (to take an $n$-vertex sample from a larger nearly optimal $\C F$-free $k$-graph) is the same as in the proofs of
Mubayi and Zhao~\cite[Proposition~1.2]{MubayiZhao07} and Lo and Markstr\"om~\cite[Proposition~1.5]{LoMarkstrom14} that the limit in~\eqref{eq:CoexLim} exists. Here, we face some new (minor) technicalities due to the fact that $\ell$-sets of zero \degree\ have to be treated differently.

We need the following auxiliary lemma which states, roughly speaking, that a $k$-graph with $n\to\infty$ vertices and positive $\ell$-\degree\ $\Omega(n^{k-\ell})$ must have $\Omega(n^k)$ edges.

\begin{lemma}\label{lm:KK} Let $0\le \ell<k$ be integers and $\gamma\ge 0$ be a real number.
If a $k$-graph $G$ with $m$ vertices and $e>0$ edges satisfies $\delta_\ell^+(G)\ge \gamma m^{k-\ell}/(k-\ell)!$ then $e\ge \gamma^{k/(k-\ell)} m^{k}/k!$.
\end{lemma}
\bpf Assume that $\gamma>0$ as otherwise there is nothing to prove. Let $\lambda>0$ be the number of $\ell$-sets covered by at least one edge of $G$. Let $x$ be the real number at least $k$ that satisfies $e={x\choose k}$, where we define 
$$
 {y\choose k}:=\frac{y(y-1)\dots (y-k+1)}{k!},\quad \mbox{for $y\in\I R$}.
 $$
 Note that the function ${y\choose k}$ is strictly increasing and continuous for $y\ge k$ so $x$ exists and is unique.
By a version of the Kruskal-Katona Theorem~\cite{Katona66,Kruskal63} that is due to Lov\'asz~\cite[Exercise 13.31(b)]{Lovasz:cpe}, we have that $\lambda\ge {x\choose \ell}$.

Thus the number of pairs $(K,L)$, where $K\in E(G)$ and $L$ is an $\ell$-subset of $K$ is  at least $\lambda\cdot \gamma m^{k-\ell}/(k-\ell)!$ on one hand and is exactly $e\cdot {k\choose \ell}$ on the other hand. Putting these two estimates together, we get that
$$
{x\choose k}{k\choose \ell}\ge {x\choose \ell}\,\frac{\gamma m^{k-\ell}}{(k-\ell)!}.
$$
 By cancelling the same (positive) factors and rearranging, we get that $(x-\ell)\dots(x-k+1)\ge \gamma m^{k-\ell}$. Now, the required inequality follows:
 $$
 k!\, e=\prod_{i=1}^k(x-i+1)\ge \left(\prod_{i=\ell+1}^k (x-i+1)\right) ^{k/(k-\ell)} \ge \gamma^{k/(k-\ell)} m^{k},
 $$
 where the first inequality can be proved by observing that, after raising it to power $k-\ell$ and cancelling identical terms, we are left with two products of $\ell(k-\ell)$ factors, with each factor at least $x-\ell+1$ on the left-hand side and  at most $x-\ell$ on the right-hand side.\epf

\bpf[Proof of Proposition~\ref{pr:LimExists}.] Let $\gamma:=\limsup_{n\to\infty} \frac{\copex_\ell(n,\C F)}{{n-\ell\choose k-\ell}}$ be the limit superior in of the stated ratios. If $\gamma=0$, then by the non-negativity of each term, the limit exists and is 0. So suppose that $\gamma>0$. Take any $\e\in (0,\gamma)$. 

Let $n$ be sufficiently large. Pick any $\C F$-free $k$-graph $G$ with $N\ge n$ vertices such that $\delta_\ell^+(G)\ge (\gamma-\e/2) {N-\ell\choose k-\ell}$.

Take a uniformly random $n$-subset of $V(G)$ and a uniformly random enumeration $v_1,\dots,v_n$ of its vertices. Equivalently, for $i=1,\dots,n$, let $v_i$ be a random element of $V(G)\setminus \{v_1,\dots,v_{i-1}\}$ with all $N-i+1$ choices being equally likely. 
 
Let $H$ be the $k$-graph on $[n]$ where a $k$-set $\{i_1,\dots,i_k\}\subseteq [n]$ is an edge if and only if $\{v_{i_1},\dots,v_{i_k}\}$ is an edge of $G$. Thus, up to relabelling of its vertices, $H$ is the $k$-graph induced in $G$ by a uniformly random set of $n$ vertices. Clearly, $H$ is $\C F$-free.

\begin{claim}\label{cl:SmallDeg} Let $g:=\big\lfloor (\gamma-\e){n-\ell\choose k-\ell}\big\rfloor$. Then the probability that there is an $\ell$-set $L\subseteq [n]$ with $0<\deg_L(H)\le g$ is less than $1/2$.
\end{claim}

\bcpf It is enough to show by the Union Bound that, for every $\ell$-set $L\subseteq [n]$, the probability over the random choices of $v_1,\dots,v_n$ that $\deg_L(H)\in [g]$ is less than $\frac12 {n\choose \ell}^{-1}$. Fix any $L\in {[n]\choose k}$. By symmetry between the vertices of $H$, we can assume for notational convenience that $L=[\ell]$. It is enough to prove the stated bound when we condition on $A:=\{v_{1},\dots,v_\ell\}$. The conditional distribution can be obtained by picking $v_{\ell+1},\dots,v_{n}$ one by one, each being uniform in the remaining subset of $V(G)\setminus A$. 
If  $\deg_A(G)=0$, then $\deg_L(H)=0$ deterministically and the stated event cannot occur.  So suppose that $\deg_A(G)>0$ and thus it is at least $\delta_\ell^+(G)\ge (\gamma-\e/2){N-\ell\choose k-\ell}$. 
For $i=0,\dots,n-\ell$, let $X_i$ be the expectation of $\deg_L(H)$ after we have exposed $v_{\ell+1},\dots,v_{\ell+i}$. In other words, $(X_0,\dots,X_{n-\ell})$ is the vertex-exposure martingale for $\deg_L(H)$ conditioned on~$A$. Note that $X_0$ is constant and its value is
\beq{eq:X0}
 X_0=\I E[X_{n-\ell}]=\deg_A(G){n-\ell\choose k-\ell}\Big/{N-\ell\choose k-\ell}\ge (\gamma-\e/2){n-\ell\choose k-\ell},
\eeq 
while $X_{n-\ell}=\deg_L(H)$. 

Let us show that $|X_{i}-X_{i-1}|\le {n-\ell-1\choose k-\ell-1}$ for every $i\in[n-\ell]$. It is enough to prove this inequality, conditioned on every choice of $v_{\ell+1},\dots,v_{\ell+i-1}$. 
For every two different choices $u$ and $u'$ for the vertex $v_{\ell+i}$ there is a natural coupling of the follow-up processes so that $|B\bigtriangleup B'|\le 2$ always holds for the current unordered sets $B,B'\subseteq V(H)$ of the selected vertices: namely, run the process for $u$ and let the process for $u'$ choose the same vertex at each step, except if we see a vertex on which the current sets $B$ and $B'$ differ then we make these two sets equal from this step onwards. Thus the final unordered sets $\{v_{\ell+1},\dots,v_n\}$ differ in at most two places and the respective {\degree}s of $L$ differ by at most ${n-\ell-1\choose k-\ell-1}$, giving the claimed inequality.

Thus Azuma's inequality (see e.g.~\cite[Theorem 7.2.1]{AlonSpencer16pm}) gives that the probability of $X_{n-\ell}\le X_0-(\e/2) {n-\ell\choose k-\ell}$ is at most
$$
\me^{-\left(\frac{\e}2{n-\ell\choose k-\ell}\right)^2\big/\left(2(n-\ell){n-\ell-1\choose k-\ell-1}^2\right)}< \me^{-\e^2 n/(9k^2)}<\frac12 {n\choose \ell}^{-1}.
$$
 Recalling that $X_{n-\ell}=\deg_L(H)$ while the constant $X_0$ satisfies~\eqref{eq:X0}, we obtain that the probability of $\deg_H(L)\le g$ is less than $\frac12 {n\choose \ell}^{-1}$, giving the claim.\ecpf

\begin{claim}\label{cl:Empty} The probability that $H$ is empty is less than $1/2$.
\end{claim}
\bcpf  Since $n$ is large, the edge density of $G$ is by Lemma~\ref{lm:KK} at least, for example, $\beta:=(\gamma-\e)^{k/(k-\ell)}>0$. Thus the expected size of $H$ is at least $\beta {n\choose k}$. Consider the  vertex-exposure martingale $(Y_0,\dots,Y_n)$ for $|E(H)|$. Similarly as before, one can show that $|Y_i-Y_{i-1}|\le {n-1\choose k-1}$ for each $i\in [n]$. Thus if $H$ has no edges then $Y_n=|E(H)|=0$ is at least $\beta {n\choose k}$ away from its mean $Y_0$ and, again by Azuma's inequality, the probability of this is at most
$$
 \me^{-\left(\beta {n\choose k}\right)^2\big/\left(2n{n-1\choose k-1}^2\right)}<\me^{-\beta^2n/(3k^2)}<\frac12,
 $$
 as desired.\ecpf

By Claims~\ref{cl:SmallDeg} and~\ref{cl:Empty} the random $\C F$-free $k$-graph $H$ on $[n]$ is non-empty and satisfies $\delta_\ell^+(H)> (\gamma-\e){n-\ell\choose k-\ell}$  with positive probability. So at least one such choice for $H$ exists and $\copex(n,\C F)>(\gamma-\e){n-\ell\choose k-\ell}$. Since $\e>0$ was arbitrary and this inequality holds for all sufficiently large $n$, we conclude that the ratio $\copex(n,\C F)/{n-\ell\choose k-\ell}$ tends to $\gamma$, finishing the proof of the proposition.\epf

\section{Positive \degree\ via hypergraph limits}

In order to state our main result of this paper (Theorem~\ref{th:main}) we need to give various definitions related to the limit theory of hypergraphs. We generally follow the notation from~\cite{Zhao15rsa}. 

For a finite set $A$ and an integer $m\ge 1$, let 
$$
 \R{A}{m}:=\{X\subseteq A\mid 0<|X|\le m\},
 $$
 consist of all non-empty subsets of $A$ with at most $m$ elements. Also, 
let 
 $$
 \r{A}:=\R{A}{|A|}=\{X\mid \emptyset\not=X\subseteq A\}
 $$ 
 denote the set of all non-empty subsets of $A$ and 
 $$
  \rl{A}:=\R{A}{|A|-1}=\{X\mid \emptyset\not=X\subsetneq A\}
  $$ 
denote the set of all \emph{proper} subsets of $A$. If $A$ is $[m]=\{1,\dots,m\}$, then we may abbreviate $\r{[m]}$ and $\rl{[m]}$ to $\rb{m}$ and $\rlb{m}$ respectively.

For a family $\C A$ of sets, let $\B x_{\C A}\in \I R^{\C A}$ denotes the vector of reals $(x_{A})_{A\in\C A}$ indexed by~$\C A$. When $\B x$ has already been specified, we use this notation as follows. If every index in $\B x$ appears in $\C A$ then by $\B x_{\C A}$ we mean an \emph{extension} of $\B x$, that is, any vector indexed by $\C A$  which coincides with $\B x$ on the common set of indices. If every element of $\C A$ appears as an index in $\B x$ then $\B x_{\C A}$ means the \emph{restriction} of $\B x$ to $\C A$, that is, the vector indexed by $\C A$ whose $A$-coordinate is the same as the $A$-coordinate of $\B x$ for every $A\in\C A$. We assume that the sets in $\C A$ come in some fixed and consistently used order, which is preserved when we pass to subfamilies. In all concrete examples for $\C A\subseteq \rb{n}$ that we give, we first order the sets increasingly by their size and then use the lexicographic order to break ties on sets of equal size.

The \emph{symmetric group} $\Sym{k}$ (consisting of all permutations of $[k]$) acts naturally on $[0,1]^{\rlb{k}}$. A function $W:[0,1]^{\rlb{k}}\to\I R$ is called \emph{symmetric} if its values do not change under the action of $\Sym{k}$ on $[0,1]^{\rlb{k}}$. For example, for $r=2$, this means that $W(x_1,x_2)=W(x_2,x_1)$ for all $(x_1,x_2)\in [0,1]^2$, and for $r=3$ this means that
 \beq{eq:SymmFor3}
 W(x_1,x_2,x_3,x_{12},x_{13},x_{23})=W(\B x^{\sigma}),\quad \mbox{for all $\B x\in [0,1]^{\rlb{3}}$ and $\sigma\in S_3$,}
\eeq
where we abbreviate $x_{\{a_1,\dots,a_k\}}$ to $x_{a_1\dots a_k}$ and denote 
$$
\B x^\sigma:=(x_{\sigma(1)},x_{\sigma(2)},x_{\sigma(3)},x_{\sigma(1)\,\sigma(2)},x_{\sigma(1)\,\sigma(3)},x_{\sigma(2)\,\sigma(3)}).
$$

A \emph{$k$-hypergraphon} is a symmetric (Lebesgue) measurable function $W:[0,1]^{\rlb{k}}\to [0,1]$. 

Elek and Szegedy~\cite[Theorem~7]{ElekSzegedy12} (see also Zhao~\cite[Theorem~1.5]{Zhao15rsa} for a different proof) showed that $k$-hypergraphons can serve as limit objects for $k$-graphs as follows. The \emph{(homomorphism) density} of a $k$-graph $F$ in $W$ is defined as
 \beq{eq:Density}
 t(F,W):=\int_{[0,1]^{\R{V(F)}{k-1}}} 
 \prod_{A\in E(F)} W(\B x_{\rl{A}})\dd\B x.
\eeq
 A $k$-graph $G$ with vertices enumerated as $v_1,\dots,v_n$ corresponds to the $k$-hypergraphon $W^G$ constructed as follows: partition $[0,1]$ into $n$ intervals $I_1,\dots, I_n$ of length $1/n$ each, and let 
\beq{eq:WG}
 W^G(\B x):=\left\{\begin{array}{ll} 1,& \mbox{if $\{v_{i_1},\dots,v_{i_k}\}\in E(G)$, where $i_j$ is the index with $I_{i_j}\ni x_j$ for $j\in [k]$,}\\
 0,&\mbox{otherwise.}
 \end{array}\right.
 \eeq
Thus $W^G$ is a $\{0,1\}$-valued function on $[0,1]^{\rlb{k}}$ that depends on the $1$-dimensional coordinates only, naturally encoding the edge set of~$G$. 
 Note that
 \beq{eq:Homs}
 t(F,W^G)=t(F,G),\quad\mbox{for every $k$-graph $F$},
\eeq
 where $t(F,G)$ is the usual \emph{(homomorphism) density} of $F$ in $G$, which is the probability that a random function $f:V(F)\to V(G)$, with all $|V(G)|^{|V(F)|}$ choices being equally likely, sends every edge of $F$ to an edge of~$G$. Call a sequence $(G_n)_{n=1}^\infty$ of $k$-graphs \emph{convergent} if, for every $k$-graph $F$, the densities $t(F,G_n)$ converge to a limit as $n\to\infty$. One of the main results of Elek and Szegedy~\cite[Theorem~7]{ElekSzegedy12} is that, for every convergent sequence $(G_n)_{i=1}^\infty$ of $k$-graphs, there is a $k$-hypergraphon $W$, called the \emph{limit} of $(G_n)_{i=1}^\infty$, such that
 \beq{eq:Convergence}
  t(F,W)=\lim_{n\to\infty} t(F,G_n),\quad \mbox{for every $k$-graph $F$}.
  \eeq

Note that, even though the right-hand side of~\eqref{eq:Convergence} involves hypergraphons depending on the 1-dimensional coordinates only (when we replace $t(F,G_n)$ by~$t(F,W^{G_n})$ using~\eqref{eq:Homs}), the resulting limit may in general depend 
on the extra (more than 1-dimensional) coordinates. Informally speaking, we may need 
to account for the limits of hypergraph constructions where the density of edges of $E(G)$ depends not only on the locations of vertices in a regularity partition of $G$ but also on the higher-dimensional cylinder structure, that is, the ``colours" of subsets of sizes between $2$ and~$k-1$. One such example is the \emph{directed cycle $3$-graph construction} $C(T)$ where one takes a quasi-random tournament $T$ and defines the edge set of $C(T)$ to consist of all triples spanning a directed cycle; this construction is known to asymptotically optimal for some extermal 3-graph problems, see e.g.~\cite{FalgasPikhurkoVaughanVolec21,GlebovKralVolec16,ReiherRodlSchacht18jems}. 
A possible corresponding $3$-hypergraphon can be defined to be $0$ except it is 1 on $\B x^{\sigma}$ for  all $\sigma\in S_3$ and $\B x$ such that $x_1<x_2<x_3$, $x_{12},x_{23}\in [0,1/2]$ and $x_{13}\in (1/2,1]$.

Accordingly, when we extend the definition of (positive) $\ell$-\degree\ to a $k$-hypergraphon $W$, we need to take into account not only 1-dimensional coordinates ($\ell$ of them) but all higher-dimensional ones (within the corresponding $\ell$-set). Formally,  for $0\le \ell<k$,  the \emph{\degree}  of $\B x\in [0,1]^{\rb{\ell}}$ in $W$ is defined as
 \beq{eq:CodegreeInW}
 \deg_W(\B x):=\int_{[0,1]^{\rlb{k}\setminus \rb{\ell}}} 
W(\B x_{\rlb{k}})\dd \B x_{\rlb{k}\setminus \rb{\ell}},
 \eeq
 that is, we take the average of $W$ over all extensions $\B x_{\rlb{k}}\in [0,1]^{\rlb{k}}$ of $\B x\in [0,1]^{\rb{\ell}}$, where the new coordinates (those indexed by ${\rlb{k}\setminus \rb{\ell}}$) are independent and uniformly distributed in~$[0,1]$.
Note that by Fubini-Tonelli's theorem (see e.g.~\cite[Theorem~2.3.2]{SteinShakarchi05ra}), the integral in~\eqref{eq:CodegreeInW} is well-defined for a.e.\ choice of $\B x\in [0,1]^{\rb{\ell}}$. For those $\B x\in [0,1]^{\rb{\ell}}$ for which the integral is undefined, we set $\deg_W(\B x):=0$ for definiteness.
 
For example, if $k=3$ and $\ell=2$, then the \degree\ of $\B x=(x_1,x_2,x_{12})$ is
 $$
 \deg_W(\B x):=\int_{[0,1]^3} W(x_1,x_2,x_3,x_{12},x_{13},x_{23})\dd x_3\dd x_{13}\dd x_{23}.
 $$
 If $W=W^G$ for a $k$-graph $G$ on $\{v_1,\dots,v_n\}$ then, with $i_j$ being the unique index such that $x_j\in I_{i_j}$ for $j\in [\ell]$ (as in~\eqref{eq:WG}), we have that 
  \beq{eq:Degs}
  \deg_{W^G}(\B x)=\left\{\begin{array}{ll} 
  \deg_G(\{v_{i_1},\dots,v_{i_{\ell}}\})\, \frac{(k-\ell)!}{n^{k-\ell}},& \mbox{if $i_1,\dots,i_\ell$ are pairwise distinct},\\
  0,&\mbox{otherwise.}\end{array}\right.
  \eeq
  Indeed, 
 every $(k-\ell)$-subset $\{v_{i_{\ell+1}},\dots,v_{i_{k}}\}$ that makes an edge of $G$ with $\{v_{i_1},\dots,v_{i_\ell}\}$ corresponds to $(k-\ell)!$ products (over all possible permutations of $\{i_{\ell+1},\dots,i_{k}\}$) of the corresponding $k-\ell$ distinct intervals $I_v$, with each product having measure~$1/n^{k-\ell}$. 
 
Call a $k$-hypergraphon $W$ \emph{non-zero} if the measure of $\B x\in [0,1]^{\rlb{k}}$ with $W(\B x)>0$ is positive. Let the \emph{minimum positive $\ell$-\degree} of a non-zero $k$-hypergraphon $W$  be defined as 
  \beq{eq:MinPCodegreeW}
   \delta_{\ell}^+(W):=\einf_{\B x\in [0,1]^{\rb{\ell}}\atop \deg_W(\B x)>0}\deg_W(\B x)=\sup_{A\subseteq [0,1]^{\rb{\ell}}\atop \mu(A)=0}  \inf\left\{\deg_W(\B x)\mid \B x\in [0,1]^{\rb{\ell}}\setminus A,\ \deg_W(\B x)>0\right\},
   \eeq
 the \emph{essential infimum} (that is, the infimum after ignoring a set $A$ of $\B x$ measure 0) of the \degree{s} $\deg_W(\B x)$ which are positive. If $W$ is zero, then we define $\delta_\ell^+(W):=0$.
 Note that, for every $k$-graph $G$, we have by~\eqref{eq:Degs} that
 \beq{eq:CodegreeWG}
 \delta_{\ell}^+(G)=\frac{(k-\ell)!}{|V(G)|^{k-\ell}}\, \delta_{\ell}^+(W^G).
 \eeq
 
For $0\le \ell<k$ and a $k$-graph family $\C F$, let $\C W^+_\ell(\C F)$ consist those $k$-hypergraphons $W$ which are the limits of some sequence of \emph{almost extremal} $k$-graphs, that is,  a sequence $(G_n)_{n=1}^\infty$ of $\C F$-free $k$-graphs such that, as $n\to\infty$, we have $|V(G_n)|\to\infty$ and $\delta_\ell^+(G_n)=(\gamma_\ell^+(\C F)+o(1)){|V(G_n)|-\ell\choose k-\ell}$. 
Also, a $k$-hypergraphon $W$ is called \emph{$\C F$-free} if $t(F,W)=0$ for every $F\in\C F$.
With this preparation, we can now formulate our main result which expresses the limit in Proposition~\ref{pr:LimExists} as the value of an optimisation problem involving $k$-hypergraphons; in fact, we give an alternative description of the set~$\C W^+_\ell(\C F)$.
 
 \begin{theorem}\label{th:main} Take any integers $0\le \ell<k$ and a  (possibly infinite) family $\C F$ of $k$-graphs. Define $\gamma:=\gamma_\ell^+(\C F)$. Then the following statements hold.
 \begin{enumerate}
 \item  The value of $\gamma$  is the supremum (in fact, maximum) of $\delta_{\ell}^+(W)$ over all $\C F$-free $k$-hypergra\-phons~$W$.
 \item A $k$-hypergraphon $W$ belongs to $\C W^+_\ell(\C F)$ if and only if it is $\C F$-free and satisfies $\delta_\ell^+(W)=\gamma$.
 \end{enumerate}
  \end{theorem}

 
Theorem~\ref{th:main} will be a direct consequence of the following two lemmas. In order to state them, we need a few more definitions. The \emph{$n$-sample} of $W$ is the distribution $\I G(n,W)$ on (vertex-labelled) $k$-graphs on $[n]$ where we sample $G\sim \I G(n,W)$ in the following two steps. First, we sample a uniform $\B x\in [0,1]^{\R{[n]}{k-1}}$ (i.e.\ each $x_A$ is uniform in $[0,1]$ and the choices over all different $A\subseteq [n]$ are mutually independent). Second, every $k$-subset $A=\{i_1,\dots,i_k\}$ of $[n]$ is included into $E(G)$ with probability $W(\B x_{\rl{A}})\in [0,1]$, with all ${n\choose k}$ choices being mutually independent. (Recall that $\B x_{\rl{A}}$ denotes the sub-vector of $\B x\in [0,1]^{\rlb{k}}$ where we take all $x_B$ with $\emptyset \not=B\subsetneq A$.) For example, if $k=3$ then $\{u,v,w\}\subseteq [n]$ is made an edge with probability
$W(x_u,x_v,x_w,x_{uv},x_{uw},x_{vw})$. One relation between $G\sim\I G(n,W)$ and the densities in $W$ is that, for every $k$-graph $F$ on $[n]$, we have
 \beq{eq:Relation}
 t(F,W)=\I P\left[ E(F)\subseteq E(G) \right],
 \eeq
 that is, 
  the $t(F,W)$ is the probability that every edge of $F$ is an edge of $G$.

Now, we are ready to state the two key lemmas and show how they imply Theorem~\ref{th:main}.

\begin{lemma}\label{lm:CodegLim} Let $k> \ell\ge 0$. Let $(G_n)_{n=1}^\infty$ be an arbitrary sequence of $k$-graphs  convergent to a $k$-hypergraphon $W$ such that $|V(G_n)|\to \infty$ as $n\to\infty$. 
Then
 \beq{eq:CodegLim}
  \delta_{\ell}^+(W)\ge \limsup_{n\to\infty} \frac{\delta_{\ell }^+(G_n)}{{|V(G_n)|-\ell\choose {k-\ell}}}.
  \eeq
  \end{lemma}

This lemma states that, informally speaking, the (normalised) minimum positive $\ell$-degree does not decrease when we pass to the limit. If $\ell\ge 1$ then we have only one-sided inequality here because there may be $o(1)$-fraction of  ``outlier'' $\ell$-tuples in $G_n$ whose \degree\ is positive but strictly smaller than $(\delta_{\ell}^+(W)+o(1)) {|V(G_n)|\choose k-\ell}$; these  $\ell$-tuples bring the positive $\ell$-\degree\ of $G_n$ down but leave no trace in the limit~$W$.

\begin{lemma}\label{lm:CodegSample}
 Let $k>\ell\ge 0$ and let $W$ be any $k$-hypergraphon. For every $\e>0$, there is $n_0$ such that for all $n\ge n_0$, if $G\sim \I G(n,W)$ then 
 the probability that 
 $|\delta_\ell^{+}(G)-\delta_\ell^{+}(W){n-\ell\choose k-\ell}|>\e{n-\ell\choose k-\ell}$ 
is at most~$\e$.
 \end{lemma}

This lemma states that minimum positive $\ell$-\degree\ of a non-zero hypergraphon is inherited within additive error $o(1)$ by a typical $n$-sample as $n\to\infty$.

\bpf[Proof of Theorem~\ref{th:main}.]%
First, suppose that $\gamma>0$. 
Take any sequence $(G_n)_{n=1}^\infty$ of almost extremal $k$-graphs.
By the standard diagonalization argument run over all (countably many) non-isomorphic $k$-graphs $F$, pass to a convergent subsequence (where $t(F,G_n)$ tends to a limit as $n\to\infty$ for each $F$). By the result of Elek and Szegedy~\cite[Theorem~7]{ElekSzegedy12} (see also~\cite[Theorem~1.5]{Zhao15rsa}) there is a $k$-hypergraphon $W$ such that $\lim_{n\to\infty} t(F,G_n)=t(F,W)$ for every $k$-graph~$F$. Of course, $t(F,W)=\lim_{n\to\infty} t(F,G_n)=0$ for every $F\in \C F$, that is, $W$ is $\C F$-free.
By Lemma~\ref{lm:CodegLim}, it holds that $\delta_{\ell}^+(W)\ge \gamma$. 

Let us show that for every $\C F$-free $k$-hypergraphon $U$, we have
 \beq{eq:U}
 \delta_{\ell}^+(U)\le \gamma.
 \eeq
 By Lemma~\ref{lm:CodegSample}, a typical $n$-vertex sample $H_n$ of $U$ for all large $n$ has minimum positive \degree\ at least $(\delta_\ell^+(U)+o(1)){n-\ell\choose k-\ell}$. (In fact, we only need this one-sided estimate from Lemma~\ref{lm:CodegSample}.) Also, for every $F\in\C F$ the probability that $H_n$ contains $F$ as a subgraph is at most~$n!\, t(F,U)=0$ by~\eqref{eq:Relation}. It follows that $(\delta_\ell^+(U)+o(1)){n-\ell\choose k-\ell}\le \copex(n,\C F)$. We conclude  that~\eqref{eq:U} holds.

By applying~\eqref{eq:U} to $W$, the limit of almost extremal $k$-graphs $G_n$, we conclude that~$\delta_\ell^+(W)=\gamma$. This proves the first part and the forward implication in the second part (also reproving Proposition~\ref{pr:LimExists}).

Let us show the converse implication in the second part. Let an $\C F$-free $k$-hypergraphon $U$ satisfy $\delta_\ell^+(U)=\gamma$. The sequence of random independent samples $(H_n)_{n=1}^\infty$, $H_n\sim \I G(n,U)$, converges to $U$ with probability $1$ by~\cite[Theorem~12]{ElekSzegedy12}. Furthermore, by Lemma~\ref{lm:CodegSample}, for every $m\ge 1$ we can find $N_m$   such that the probability of $\delta_\ell^+(H_{N_m})< (\gamma-1/m){N_m-\ell\choose k-\ell}$ is at most~$2^{-m-1}$. Clearly, $H_{N_m}$ is $\C F$-free with probability 1. Thus, with probability at least $1-\sum_{m=1}^\infty 2^{-m-1}>0$, $(H_{N_m})_{m=1}^\infty$ is a sequence of almost extremal $k$-graphs convergent to $U$, that is, $U\in\C W^+_\ell(\C F)$, proving the second part.

Finally, suppose that $\gamma=0$. The only non-trivial claim that we have to establish is that no $\C F$-free $k$-hypergraph $U$ can satisfy $\delta_\ell^+(U)>0$ and this follows as above by taking random samples from $U$ and applying Lemma~\ref{lm:CodegSample}.\epf

\hide{
Let us remark that it seems challenging to prove the first part of Theorem~\ref{th:main} without recoursing to hypergraphons. The difficulty of the naive approach of 
taking a sequence $(G_m)$ of $k$-graphs that attain the lim-sup in~\eqref{eq:Gamma+}
and applying Lemma~\ref{lm:CodegSample} to (the hypergraphon of) $G_m$ for large $m=m(n)$ is that the orders of the $k$-graphs $G_m$ can in principle be very far spaced from each other and then the straightforward Union Bound over all $m^\ell$-subsets $L$ of $[n]$ of size $\ell$ does not work directly, since the probability that a particular choice $L$ fails can in general be bounded by a function of $n$ only. 
}

\section{Proofs of Lemma~\ref{lm:CodegLim} and~\ref{lm:CodegSample}}

In order to present the proofs of Lemma~\ref{lm:CodegLim} and~\ref{lm:CodegSample} we need some further notation.

The following definition of a partially vertex-labelled hypergraph will suffice for the purposes of this paper. Namely, for $k>\ell\ge 0$,
an \emph{$\ell$-labelled} $k$-graph is a triple $F=(V,E,\ell)$ where $(V,E)$ is a $k$-graph and $V\supseteq [\ell]$. We view $1,\dots,\ell\in V$ as labelled vertices and call them the \emph{roots}. 

For a $k$-hypergraphon $W$, a vector $\B x\in [0,1]^{\rb{\ell}}$ and an $\ell$-labelled $k$-graph $F=(V,E,\ell)$, the \emph{($\B x$-rooted) density} of $F$ in $W$ is
 \beq{eq:RootedDensity}
 t_{\B x}(F,W):=\int_{[0,1]^{\R{V}{k-1}\setminus \rb{\ell}}} \prod_{A\in E} W(\B x_{\rl{A}})\dd \B x_{\R{V}{k-1}\setminus \rb{\ell}}.
 \eeq
 By Fubini-Tonelli's theorem, this is defined for a.e.\ $\B x\in [0,1]^{\rb{\ell}}$; for all other $\B x$ we set $t_{\B x}(F,W):=0$ for definiteness.
 For example, we have that  the definition in~\eqref{eq:RootedDensity} for 
 $$
 \Edge{k}{\ell}:=([k],\{[k]\},\ell),
 $$ 
 the $\ell$-labelled single $k$-edge, becomes exactly the one in~\eqref{eq:CodegreeInW}, and thus
 $$
 t_{\B x}(\Edge{k}{\ell},W)=\deg_W(\B x),\quad\mbox{for every $\B x\in [0,1]^{\rb{\ell}}$},
 $$
The densities of an $\ell$-labelled $k$-graph $F=(V,E,\ell)$ and its unlabelled version $\eval{F}:=(V,E)$, where we just forget the labelling, satisfy by Fubini-Tonelli's theorem the following relation:
 \beq{eq:AvOfRootedD}
  t(\eval{F},W)=\int_{[0,1]^{\rb{\ell}}} t_{\B x}(F,W)\dd\B x.
  \eeq
  This can be informally interpreted as that the density of the unlabelled $k$-graph $\eval{F}$ is the average over the uniform choice of an ``$\ell$-tuple'' $\B x\in [0,1]^{\rb{\ell}}$ of the $\B x$-rooted density of~$F$.

 The \emph{product} $FF'$ of any
 two $\ell$-labelled $k$-graphs $F=(V,E,\ell)$ and $F'=(V',E',\ell)$ is obtained by replacing $F'$ by an isomorphic $\ell$-labelled $k$-graph with $V\cap V'=[\ell]$ (that is, making $F$ and $F'$ vertex-disjoint except for the roots) and then taking the union of the vertex and edge sets: 
 $$
  FF':=(V\cup V',E\cup E,\ell).
  $$
   The name ``product'' comes from the relation 
 \beq{eq:Product}
 t_{\B x}(FF',W)=t_{\B x}(F,W)\,t_{\B x}(F',W),\quad\mbox{for every $\B x\in [0,1]^{\rb{\ell}}$},
 \eeq
 which holds since the integral for $t_{\B x}(FF',W)$ can be written by Fubini-Tonelli's theorem as the product of two integrals by partioning its variables into two groups: namely $x_A$ with $A\cap (V\setminus [\ell])\not=\emptyset$ and $x_{A'}$ with $A'\cap (V'\setminus [\ell])\not=\emptyset$. (Recall that $\ell<r$ so the roots do not span any edges.) It follows that, with $\Edge{k}{\ell}^m$ denoting the $m$-fold product of the $\ell$-labelled single $k$-edge with itself, we have
 \beq{eq:Power}
 t_{\B x}(\Edge{k}{\ell}^m,W)=\left(t_{\B x}(\Edge{k}{\ell},W)\right)^m,\quad\mbox{for every $\B x\in [0,1]^{\rb{\ell}}$}.
  \eeq


Now we are ready to prove the lemmas.
 
\bpf[Proof of Lemma~\ref{lm:CodegLim}.] Suppose that a sequence $(G_n)_{n=1}^\infty$ convergent to some $W$ gives a counterexample to the lemma.  By passing to a subsequence, we can assume that the ratios $\delta_{\ell }^+(G_n)/{|V(G_n)|-\ell\choose {k-\ell}}$ converge to some $\delta$, which is strictly larger that $\delta_\ell^+(W)$. By Lemma~\ref{lm:KK}, the edge density of $G_n$ is at least $\delta^{k/(k-\ell)}+o(1)>0$ as $n\to\infty$, so their limit $W$ is non-zero. Thus the set 
$$
X:=\{\B x\in [0,1]^{\rb{\ell}}\mid \deg_W(\B x)\in (0,\delta)\}
$$ has positive measure.

By the countable additivity of measure, 
there is $\e>0$ such that the set 
$$
 X_{\e}:=\left\{\B x\in [0,1]^{\rb{\ell}}\mid \deg_W(\B x)\in (\e,\delta-\e)\right\},
 $$
 has measure at least $\e$. Indeed, $X$ is the countable union $\cup_{m=1}^\infty X_{1/m}$, so $X_{1/m}$ has positive measure for some $m$ and we can take $\e:=\min\{1/m,\mu(X_{1/m})\}$. 
 
\begin{figure}
\begin{center}
\includegraphics[scale=1]{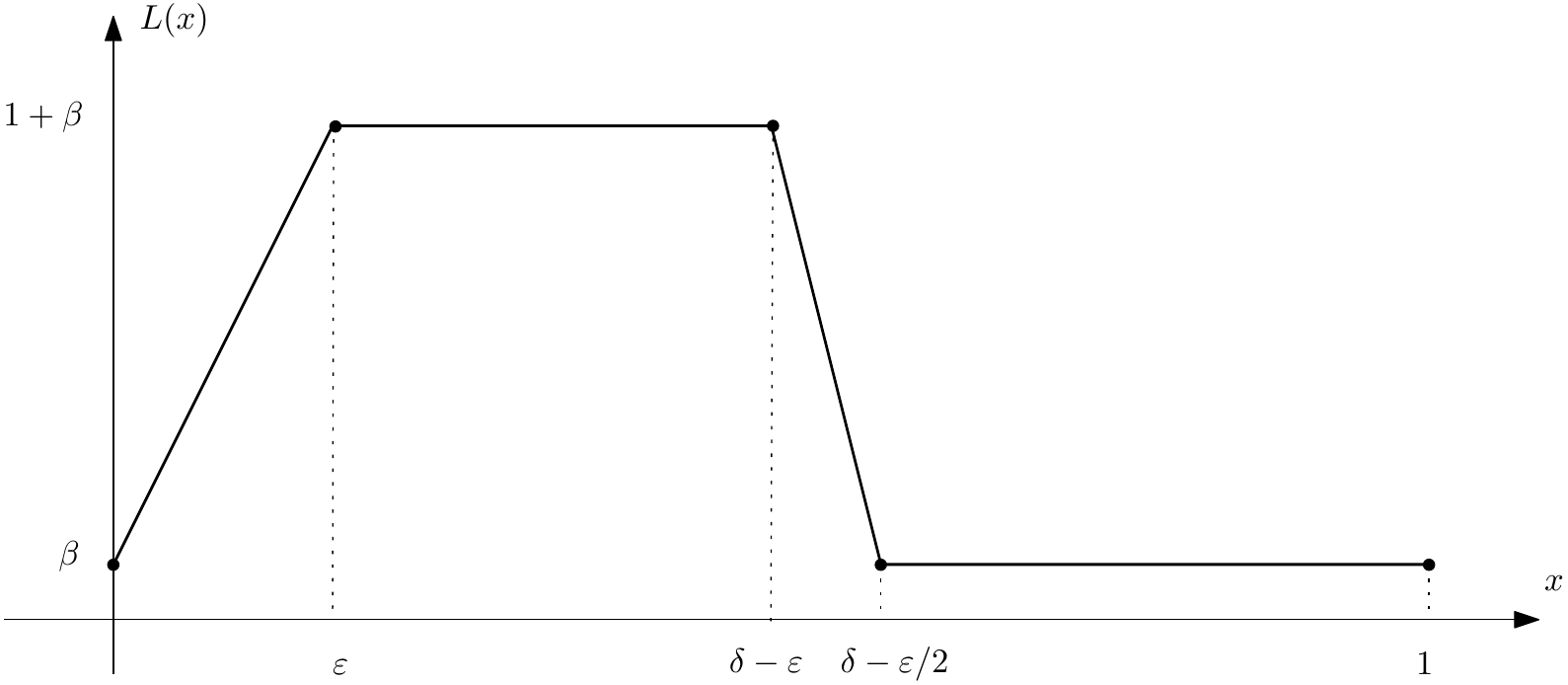}
\end{center}
\caption{The graph of the function $L$}\label{fi:1}
\end{figure}

 Since $X_{\e}\not=\emptyset$, we have $\e< \delta/2$.
Fix any $\beta\in (0,\e/2)$. Let $L(x):[0,1]\to \I R$ be the piecewise linear function whose graph in $[0,1]\times\I R$ consists of the linear segments connecting the points $(0,\beta)$, $(\e,1+\beta)$, $(\delta-\e,1+\beta)$, $(\delta-\e/2,\beta)$, and $(1,\beta)$ in this order, see Figure~\ref{fi:1}.
In particular, $L$ is constant $1+\beta$ on $[\e,\delta-\e]$ and constant $\beta$ on $[\delta-\e/2,1]$. Informally, the ``penalty'' function $L$ penalises values strictly between $0$ and $\delta$.

Since the function $L$ is continuous, the Stone-Veierstrass Theorem gives a polynomial $p(x)=\sum_{i=1}^D a_i x^i$ such that $|p(x)-L(x)|\le \beta$ for every $x\in [0,1]$. This polynomial clearly has the following properties:
 \begin{eqnarray}
  p(x)\ge  0, && \mbox{for all $x\in [0,1]$},\label{eq:PNonNeg}\\
  p(x)\le 2\beta, && \mbox{for all $x\in \{0\}\cup [\delta-\e/2,1]$},\label{eq:PSmall}\\
  p(x)\ge 1,&&  \mbox{for all $x\in[\e,\delta-\e]$}.\label{eq:PLarge}
  \end{eqnarray}

For a $k$-hypergraphon $U$, define $Q_U:=\int_{[0,1]^{\rb{\ell}}} q_U(\B x)\dd \B x$ to be the average of 
$$
 q_U(\B x):=p\left(t_{\B x}(\Edge{k}{\ell}, U)\right)
 $$
 over uniform $\B x\in [0,1]^{\rb{\ell}}$. Note by~\eqref{eq:Power} that
   $$
   q_U(\B x)=\sum_{i=0}^D a_i \left(t_{\B x}(\Edge{k}{\ell},U)\right)^i = a_0+\sum_{i=1}^D a_i\, t_{\B x}(\Edge{k}{\ell}^i,U),\quad\mbox{for all $\B x\in [0,1]^{\rb{\ell}}$},
   $$
   Thus, by the linearity of integral and by~\eqref{eq:AvOfRootedD},
   \begin{eqnarray*}
   Q_U&=&a_0+\sum_{i=0}^D a_i  \int_{[0,1]^{\rb{\ell}}} t_{\B x}(\Edge{k}{\ell}^i,U)\dd\B x = a_0+\sum_{i=1}^D a_i\,   t(\eval{\Edge{k}{\ell}^i},U).
   \end{eqnarray*}
   
   Let $U_n:=W^{G_n}$ be the hypergraphon of $G_n$.
   As $G_n$ converges to $W$, we have that, for every $i\in[D]$, the (unlabelled) $k$-graph density $t(\eval{\Edge{k}{\ell}^i},U_n)$ converges to $t(\eval{\Edge{k}{\ell}^i},W)$ and thus
   \beq{eq:QGn}
   \lim_{n\to\infty} Q_{U_n}=Q_W.
   \eeq

If we take $n\to\infty$ and evaluate $Q_{U_n}$ then, with $m:=|V(G_n)|$, the outer integral becomes the average value of the polynomial $p$ evaluated at the (obviously defined) rooted density $t_{v_1,\dots,v_\ell}(\Edge{k}{\ell},G_n)$ where $v_1,\dots,v_\ell$ are independent uniformly chosen vertices of~$G_n$.  For each of these evaluations of $p$, its argument is either $0$ (if some two $v_i$'s coincide or the $\ell$-set $\{v_1,\dots,v_\ell\}$ is not covered by any edge of $G_n$) or at least $\delta_{\ell}^+(G_n)\cdot (k-\ell)!/m^{k-\ell}$ (in all other cases). 
Thus, if $n$ is large enough, then by $\delta_{\ell}^+(G_n)\cdot (k-\ell)!/m^{k-\ell}=\delta+o(1)>\delta-\e/2$ and by~\eqref{eq:PSmall} each computed value of $p$ is at most $2\beta$. Thus we have that $Q_{{U_n}}\le 2\beta$ for all large~$n$.

On the other hand, since $q_W(\B x)\ge 0$ for every $\B x$ by~\eqref{eq:PNonNeg}, we have by~\eqref{eq:PLarge} that
 $$
 Q_W\ge \int_{X_\e} q_W(\B x)\dd \B x\ge 1\cdot \mu(X_\e)\ge \e.
 $$  
 Since $2\beta<\e$, $Q_W$ cannot be be the limit of $Q_{U_n}$, a contradiction to~\eqref{eq:QGn} proving Lemma~\ref{lm:CodegLim}.\epf


\bpf[Proof of Lemma~\ref{lm:CodegSample}.] Let $\delta:=\delta_\ell^+(W)$, $n\to\infty$ and let $G\sim \I G(n,W)$. We have to show that $\delta_\ell^+(G)$ is unlikely to be far from $\delta{n-\ell\choose k-\ell}$.
 Assume that $W$ is non-zero as otherwise $G$ is empty with probability 1 and the lemma trivially holds. 

Consider the vertex exposure martingale $(Y_0,\dots,Y_n)$ with $Y_n=|E(G)|$, where for $i=1,\dots,n$ we expose all $x_A$'s with $\max A= i$ as well as all edges of $G$ whose maximal element is~$i$. We have that $|Y_i-Y_{i-1}|\le {n-1\choose k-1}$ for every $i\in [n]$ because if we change our choices at Step~$i$ this will affect only edges of $G$ containing $i$. (Note that, unlike in the proof of Proposition~\ref{pr:LimExists}, we do not have to worry about the measure-0 event that some different $x_i$'s coincide.) Also, 
$$
 Y_0=\I E(Y_n)=t(\Edge{k}{0},W){n\choose k},
 $$
  which is $\Omega(n^{k})$ since $W$ is non-zero.
Azuma's inequality gives that the probability that $G$ spans no edges (i.e.\ $Y_n=0$) is $\me^{-\Omega(n)}<\e/3$. 


Next, let us show that, for every fixed $\ell$-set $L\subseteq [n]$, the probability that its \degree\ $\deg_L(G)$ in $G\sim \I G(n,W)$ is positive but less than $(\delta-\e){n-\ell\choose k-\ell}$ is at most $(\e/3){n\choose \ell}^{-1}$. By symmetry, assume that $L=[\ell]$.%
\hide{Let 
 $$
 G_L:=\{X\subseteq [n]\setminus L\mid X\cup L\in E(G)\},
 $$
 be the link of $L$ in the random $k$-graph $G$.} 
 Take any $\B x\in [0,1]^{\rb{\ell}}$. By ignoring a set of $\B x$ of measure 0, we have that $\deg_W(\B x)$ is 0 or at least $\delta$. In the former case, $\deg_L(G)=0$ with probability $1$. So suppose that $\deg_W(\B x)\ge \delta$.
Given $\B x$, consider the natural vertex exposure martingale $(X_0,\dots,X_{n-\ell})$ for $\deg_L(G)$. It holds that $|X_i-X_{i-i}|\le {n-\ell-1\choose k-\ell-1}$ for each $i\in [n-\ell]$. Since 
$$
 X_0=\I E[X_{n-\ell}]=\deg_W(\B x){n-\ell\choose k-\ell}\ge \delta{n-\ell\choose k-\ell},
 $$
  Azuma's inequality gives that the probability of $X_{n-\ell}<X_0-\e {n-\ell\choose ^{k-\ell}}$ is $\me^{-\Omega(n)}< (\e/3){n\choose \ell}^{-1}$, as claimed. 

The Union Bound over all ${n\choose \ell}$ choices of an $\ell$-set $L\subseteq [m]$ shows that the probability that $G$ is non-empty and $\delta_\ell^+(G)< (\delta -\e){n-\ell\choose ^{k-\ell}}$ is at most~$2\e/3$. 

Finally, it remains to upper bound the probability that $\delta_\ell^+(G)$ is too large. Since this part of the lemma is not used anywhere else in the paper, we will be rather brief. By the definition of $\delta=\delta_\ell^+(W)$ and since $W$ is non-zero, the set 
 $$
 Y:=\{\B x\in [0,1]^{\rb{l}}\mid 0<\deg_W(\B x)\le \delta+\e/2\}
 $$
 has positive measure. By decreasing $\e>0$ if necessary, assume that 
 $$
  Z:=\{\B x\in [0,1]^{\rb{l}}\mid \e<\deg_W(\B x)\le \delta+\e/2\}
  $$
  has positive measure.
 Take a uniform $\B x\in [0,1]^{\R{[n]}{k-1}}$. The expected number of $\ell$-sets $L\subseteq [n]$ such that $\B x_{\r{L}}\in X$ is $\Omega(n^\ell)$. By Azuma's inequality, the probability that no $L$ satisfies this is $\me^{-\Omega(n)}<\e/6$. Furthermore, if we take any $L$ with $\B x_{\r{L}}\in Z$ and condition on $\B x_{\r{L}}$ then Azuma's equality, when applied to the martingale where we expose for each vertex $i\in [n]\setminus L$ all information up to $i$, shows that the probability of $\deg_G(L)$ being at least $(\e/2){n-\ell\choose k-\ell}$ away from its expected value $\deg_W(\B x){n-\ell\choose k-\ell}$ is $\me^{-\Omega(n)}<\e/6$. Thus with probability at least $1-\e/3$ there is an $\ell$-set $L\subseteq [n]$ with $0<\deg_G(L)\le (\delta+\e){n-\ell\choose k-\ell}$, which implies that $\delta_\ell^+(G)\le (\delta+\e){n-\ell\choose k-\ell}$.

By putting all together, we conclude that, with probability at least $1-\e$, the random $n$-sample $G\sim \I G(n,W)$  has positive $\ell$-\degree\ $(\delta\pm\e){n-\ell\choose ^{k-\ell}}$, as desired.\epf

\section{Concluding remarks}

One can also consider extremal limits for the $\ell$-degree Tur\'an problem. Namely, let
$\C W_\ell(\C F)$ consist of all $k$-hypergraphons $W$ such that there is a sequence $(G_n)_{n=1}^\infty$ of $\C F$-free $k$-graphs convergent to $W$ such that $|V(G_n)|\to\infty$ and the \emph{mininum $\ell$-degree}
$$\delta_\ell(G_n):=\min\{\deg_{G_n}(L)\mid L\subseteq [n],\ |L|=\ell\}
$$
 is $(\gamma_\ell(\C F)+o(1)){|V(G_n)|-\ell\choose k-\ell}$ as $n\to\infty$, where $\gamma_\ell(\C F)$ is defined in~\eqref{eq:CoexLim}. Also, let the \emph{minimum $\ell$-degree} of a $k$-hypergraphon $W$ be
 $$
 \delta_\ell(W):=\einf_{\B x\in [0,1]^{\rb{\ell}}} \deg_W(\B x).
 $$
 Our proof of Theorem~\ref{th:main} can be easily adapted to to produce the following result (in fact, the proof is simpler since we do not have to treat $\ell$-sets of zero degree in a special way).

\begin{theorem} Let $k>\ell\ge 0$ be integers and let $\C F$ be a $k$-graph family. Then $\gamma_\ell(\C F)$ is the maximum of $\delta_\ell(W)$ over all $\C F$-free $k$-hypergraphons $W$. Moreover,  $W\in\C W_\ell(\C F)$ if and only if $W$ is $\C F$-free and $\delta_\ell(W)=\gamma_\ell(\C F)$.
 \end{theorem}
 
\hide{ 
In general, even if $\gamma_\ell(\C F)=0$, the set $\C W_\ell(\C F)$ may consists of non-zero graphons. For example, if $\C F$ consists of the $3$-graph on $[5]$ with edges $\{1,2,3\}$, $\{1,2,4\}$ and $\{3,4,5\}$ then it is easy to see that $\gamma_\ell(\C F)=0$ but uniform blowups of a single edge are $\C F$-free and have edge density $2/9$. However, the situation with the positive degree Tur\'an problem is unclear:

\begin{problem} Does  $\gamma_\ell^+(\C F)=0$ always imply that every $\C F$-free graph $W$ is zero?\end{problem}
}



\end{document}

\section*{Macros}

$\coex(n,\C F)$ \verb$\coex(n,\C F)$

$\copex(n,\C F)$ \verb$\copex(n,\C F)$

$\r{A}$ \verb$\r{A}$

$\R{A}{m}$ \verb$\R{A}{m}$

$\rb{k}$ \verb$\rb{k}$

$\rl{A}$ \verb$\rl{A}$

$\rlb{k}$ \verb$\rlb{k}$

$\Sym{m}$ \verb$\Sym{m}$

$\deg_W(\B x)$ \verb$\deg_W(\B x)$

$\essinf$ \verb$\essinf$

$\Edge{k}{m}$ \verb$\Edge{k}{m}$

$\eval{F}$ \verb$\eval{F}$

\end{document}